\documentclass[12pt]{amsart}
\usepackage{amssymb}
\usepackage{amsfonts}
\usepackage{amsmath}
\usepackage{array}
\usepackage{rotating}
\usepackage[utf8]{inputenc}
\usepackage{enumitem}
\usepackage{placeins}
\usepackage{pspicture}
\usepackage{epsf}

\usepackage{xcolor}

\usepackage{tikz}

\usepackage[matrix,ps,xdvi]{xypic} 
\newdimen\vcadre\vcadre=0.1cm 
\newdimen\hcadre\hcadre=0.1cm 

\def\arx#1[#2]{\ifcase#1 \relax \or%
  \ar @{-}[#2]  \or%
  \ar @2{-}[#2] \or%
  \ar @{--}[#2] \or%
  \ar @2{.}[#2] \or%
  \ar @{~}[#2]  \fi}

\def\ST{{\rm ST}}
\def\PST{{\rm PST}}

\def\AA{{\mathbb A}}

\def\<{\langle}
\def\>{\rangle}

\def\Z{\operatorname{\mathbb Z}}

\def\SG{{\mathfrak S}}

\usepackage{placeins}
\usepackage{pspicture}

\def\PF{{\rm PF}}

\def\Tabvrule{\vrule width-0.4pt}       
\def\Tabhrule{\hrule \hrule height-0.4pt} 
\def\Tabstrut{\vrule height2.2ex 
                     depth0.8ex  
                     width0ex    
\relax}

\def\PasCase#1{\omit%
            $\vcenter{\hbox {\vbox to 0.4pt{}}
               \hbox{\makebox[3ex]{\Tabstrut$#1$}}}%
               \Tabvrule$}
\def\PasCasePoint{\PasCase{\cdot}}
\def\DessinCarre#1{%
    \vcenter{\hbox{}\hrule
             \hbox{\vrule\makebox[3ex]{\Tabstrut$#1$}\vrule}\Tabhrule}%
             \Tabvrule}
\def\GenRuban#1{\vcenter{\halign{&$\DessinCarre{##}$\cr#1}}\egroup}

\def\sTabvrule{\vrule width-0.4pt}
\def\sTabhrule{\hrule \hrule height-0.4pt}
\def\sTabstrut{\vrule height1.6ex depth0.6ex width0ex \relax}
\def\sDessinCarre#1{%
    \vcenter{\hbox{}\hrule
             \hbox{\vrule\makebox[2.3ex]%
                  {\sTabstrut$\scriptstyle#1$}\vrule}\sTabhrule}%
             \sTabvrule}
\def\sGenRuban#1{\vcenter{\halign{&$\sDessinCarre{##}$\cr#1}}\egroup}

\def\ruban{%
  \bgroup
  \let\ =\omit
  \let\\=\cr
  \let\x=\times
  \let\.=\PasCasePoint
  \offinterlineskip
  \GenRuban}

\def\sruban{%
  \bgroup
  \let\ =\omit
  \let\x=\times
  \let\\=\cr
  \offinterlineskip
  \sGenRuban}

\newdimen\vcadre\vcadre=0.2cm 
\newdimen\hcadre\hcadre=0.2cm 

\def\cerp#1#2{\put(#1,#2){\circle*{0.7}}}
\def\cerg#1#2{\put(#1,#2){\circle*{1}}}

\def\arbtgb{\begin{picture}(3,6)\cerg25\cerp13\cerp33\cerp11
 \put(2,5){\Line(-1,-2)}
 \put(2,5){\Line( 1,-2)}
 \put(1,3){\Line( 0,-2)}
\end{picture}}

\def\cerp#1#2{\put(#1,#2){\circle*{0.7}}}
\def\cerg#1#2{\put(#1,#2){\circle*{1}}}

\title[]%
{The noncommutative geode}

\author[J.-C.~Novelli and J.-Y.~Thibon]%
{Jean-Christophe Novelli and Jean-Yves Thibon}

\address[] {Laboratoire d'Informatique Gaspard Monge, Universit\'e Gustave
Eiffel, CNRS, ENPC, ESIEE-Paris, \\
5 Boulevard Descartes \\Champs-sur-Marne \\77454 Marne-la-Vall\'ee cedex 2 \\
FRANCE}
\email[Jean-Christophe Novelli]{jean-christophe.novelli@univ-eiffel.fr}
\email[Jean-Yves Thibon]{jean-yves.thibon@univ-eiffel.fr} 
\date{}



\begin{document}

\begin{abstract}
We investigate the geode and some of its generalizations from the point of view on noncommutative symmetric functions.
\end{abstract}

\maketitle
\section{Introduction}
Let $f(t)$ be a formal power series with nonzero constant term.
By the \emph{Lagrange series}, we shall mean the (unique) formal power series (regarded as a function of $f$)
\begin{equation}
g(t)=\sum_{n\ge 0}g_nt^n
\end{equation} 
solving the functional equation
\begin{equation}\label{eq:lag}
g(t) = f(tg(t)) =
 \sum_{ n\ge 1}f_nt^ng(t)^n\ \text{where}\ f(t)=\sum_{n\ge 0}f_nt^n,\ f_0\not=0.
\end{equation}
In \cite{WR}, Wildberger and Rubine obtain a curious property of the Lagrange series (for $f_0=1$):
\begin{equation}
 \gamma(t):= \frac{g(t)-1}{f(t)-1}
\end{equation}
is a formal power series whose coefficients are polynomials in the $f_i$ with
nonnegative integer coefficients, which they call the \emph{geode}.

The aim of this note is to discuss the geode in the context of our previous works \cite{NTpark, NTLag,NTLag2,MNT,NTDup,JVNT} on the
symmetric and noncommutative versions of the Lagrange series.

\section{The symmetric and noncommutative Lagrange series}

When the $f_n$ are the complete homogenous symmetric functions $f_n=h_n(X)$,
we obtain the \emph{symmetric Lagrange series}. Then $g_n(X)$ is the Frobenius
characteristic of the permutation representation of the symmetric group
$\SG_n$ on the set $\PF_n$ of parking functions of length $n$
\cite{Hai1,NTLag}. The associated multiplicative basis of symmetric functions
$g_\lambda$ is also involved in Macdonald's formula for the top connection
coefficients  of $\SG_n$ \cite{Mcd,GJ1,GJ2,NTLag2} and in the description of
the reduced incidence Hopf algebra of the lattices of noncrossing partitions
\cite{Einz,EH,NTLag2}.

If we assume that the $f_n$ are noncommuting indeterminates, we obtain the
\emph{noncommutative Lagrange series}, which might as well be called the
Łukasiewicz series. Indeed, solving \eqref{eq:lag} recursively, we see that
$g_n$ becomes the sum of all Łukasiewicz words of length $n+1$ and sum
$n$~\cite{Ran}: writing
for short $f_{i_1i_2\cdots}$ for $f_{i_1}f_{i_2}\cdots$,
\begin{equation}
 g_0 = f_0,\
 g_1 = f_{10},\
 g_2 = f_{200}+f_{110},\
 g_3 = f_{3000}+f_{2100}+f_{2010}+f_{1200}+f_{1110},
 \ldots
\end{equation}
{\it i.e.}, the Polish codes for plane rooted trees on $n+1$ vertices (obtained
by reading the arities of the nodes in prefix order). And indeed, for $t=1$,
\eqref{eq:lag} reads
\begin{equation}
 g = f_0 + f_1g +f_2g^2 + f_3g^3+\cdots
\end{equation}
which is just the Łukasiewicz grammar
\begin{equation}
	G\rightarrow 0|1G|2GG|3GGG|\cdots
\end{equation}
encoded as a formal series.

These words also encode
in a natural way various Catalan sets. Setting $f_i=a^ib$, we obtain Dyck
words with an extra $b$ at the end. Seeing $f_{i_1i_2\cdots i_r}$ as
encoding the nondecreasing word $1^{i_1}2^{i_2}\cdots r^{i_r}$, we obtain a
nondecreasing parking function, which can itself be decoded as a noncrossing
partition, whose blocks are encoded by their minimal elements repeated as many
times as the lengths of the blocks.

\medskip
{\footnotesize
For example, the  word $f_{2100}$ encodes the plane tree $\arbtgb\ $, the
Dyck word $aababb\cdot b$, the nondecreasing parking function $112$ and the
noncrossing partition $13|2$.
}

\medskip
The existence of the geode is immediately apparent on the Łukasiewicz series.
Define the maps $d_k$ on trees as: given a tree $T$, if the last nonzero
value in the Łukasiewicz word associated with $T$ is not $k$, then send $T$ to
zero; otherwise, this last nonzero value corresponds to a corolla in $T$ that
is the last one in prefix order and then $d_k$ send $T$ to the tree
obtained by replacing this corolla with a leaf.

Applying $d_k$ to $g_{n+k}$ gives a multiset of trees $\Gamma_n$ in
$g_n$ that is independent of $k$: indeed, given a tree of size $n$, the number
of ways to add a corolla as the last one in prefix order is independent of the
size of the corolla.

For example, with $k=3$ and $n=5$, the second corolla with a $3$ as root is
the last one in prefix order so it is removed by application of $d_3$ and
yields $31010$.

\begin{equation}
{\xymatrix@C=2mm@R=2mm{
 *{} & {3}\ar@{-}[dl]\ar@{-}[d]\ar@{-}[dr] \\
 {1}\ar@{-}[d] & {1}\ar@{-}[d] & {0} & *{} & {} \\
 {0} & {3}\ar@{-}[dl]\ar@{-}[d]\ar@{-}[dr] & {} & *{} & {} \\
 {0} & {0} & {0} & *{} & {} \\
}}
\qquad
\to
\qquad
{\xymatrix@C=2mm@R=2mm{
 *{} & {3}\ar@{-}[dl]\ar@{-}[d]\ar@{-}[dr] \\
 {1}\ar@{-}[d] & {1}\ar@{-}[d] & {0} & *{} & {} \\
 {0} & {0} & {} & *{} & {} \\
}}
\end{equation}

Eliminating the zeros by setting $f_0=1$, the sum $\gamma_n$ of their
codes satisfies
\begin{equation}
  g_n = f_n+\gamma_1 f_ {n-1}+\gamma_2 f_{n-2}+\cdots+ \gamma_{n-1}f_1,
\end{equation}
which is the (noncommutative) equation of the geode.

From now on, we shall assume that $f_n=S_n$, the noncommutative complete
symmetric function so that $g$ becomes the  noncommutative symmetric Lagrange
series:
\begin{equation}
\label{g0123}
 g_0 = 1,\ 
 g_1 = S_1,\
 g_2 = S_2+S^{11},\
 g_3 = S_3+2S^{21}+S^{12}+S^{111},\dots
\end{equation}
Then $g_n(A)$ is the noncommutative Frobenius characteristic of the natural
representation of the $0$-Hecke algebra $H_n(0)$ on parking
functions~\cite{NTLag}.

\section{The noncommutative symmetric geode}

To compute the noncommutative symmetric geode, we can make use of the operators $S_n^{-1}$, defined by
\begin{equation}
	S^{i_1\ldots i_r}S_n^{-1} = 
	\begin{cases}
		S^{i_1\ldots i_{r-1}}&\text{if $i_r=n$},\\
		     0&\text{otherwise}.
	\end{cases}
\end{equation}
and since the above argument shows that $\gamma_n = g_{n+k}S_k^{-1}$, we can
compute $\gamma$ by just applying $S_1^{-1}$ to $g$, yielding 
\begin{equation}
\label{gamma0123}
 \gamma_0 = 1,\ 
 \gamma_1 = S_1,\
 \gamma_2 = 2S_2+S^{11},\
 \gamma_3 = 3S_3+3S^{21}+2S^{12}+S^{111},\dots
\end{equation}

The coefficient of $S^I$ in $g_n$ is equal to the number of ordered trees on
$n + 1$ vertices whose sequence of non-zero arities in prefix order is $I$. To
compute the coefficient of $S^I$ in $\gamma_n$, one has to add for each of
these trees the number of final zeros of their code. For example, $3000$
yields $3S^3$, $2100$ and $2010$ respectively yield $2S^{21}$ and $S^{21}$,
$1200$ gives the term $2S^{12}$, and $1110$ the term $S^{111}$.
This description is equivalent to~\cite[Theorem 1]{Gos}.

Alternatively, the number of trailing zeros of the code is also the number of
possibilties of shifting the code to the right so as to obtain the evaluation
of a word on $[n]$. For example, $3000$ can be shifted to $0300$ and $0030$,
which correspond to the words $111,222,333$. Similarly, the shifts of $2100$
yield the words $112,223$, $2010$ gives only $113$, $1200$ gives $122,233$ and
$1110$ only $123$.

Clearly, $\gamma_n$ is the noncommutative Frobenius characteristic of the
representation of $H_n(0)$ on rearrangements of these words, which can be characterized as all
shifts of nondecreasing parking functions of length $n$ restricted to
the alphabet~$[n]$.

The number of (nonempty) such words (sum of the coefficients of $\gamma_n)$ is
Sequence~A071724:
$$1,1,3,9,28,90,297,1001,\dots$$
as one can check with the example $n=3$ above. The generating function of these numbers is obtained by specializing $g$ to the alphabet $A=x$, so that $S_n\mapsto x^n$ and $g$ becomes
\begin{equation}
g(x)= C(x)=\frac{1-\sqrt{1-4x}}{2x},
\end{equation}
the generating function of Catalan numbers, and $\gamma$ becomes
\begin{equation}
\gamma(x)=\frac{(C(x)-1)(1-x)}{x}.
\end{equation}

\section{Expansion on other bases}

The coefficients of $g$ in the bases $R_I$ and $\Lambda^I$ also have interesting combinatorial
interpretations, which have been described in \cite{NTLag}. Taking into account the properties
\begin{equation}
R_{i_1\ldots i_r}S_n^{-1} = 
	\begin{cases}
		R_{i_1\ldots i_{r-1}}&\text{if $i_r=n$},\\
		     0&\text{otherwise},
	\end{cases}
\end{equation}
and
\begin{equation}
\Lambda^{i_1\ldots i_r}S_1^{-1} = \Lambda^{i_1\ldots i_{r}-1},
\end{equation}
we obtain the description of the coefficients of $\gamma$ on these bases.

The coefficient of $R_I$ in $\gamma$ is equal to the number of parking
quasi-ribbons of shape~$I1$. For example, in
\begin{equation}
\gamma_3= 9R_3 +4R_{21} + 3R_{12} + R_{111},
\end{equation}
the coefficients count
\begin{align}
R_3:& 111|2, 111|3, 111|4, 112|3, 112|4, 1213|4, 122|3, 122|4, 123|4\\
R_{21}: &11|2|3, 11|2|4, 11|3|4, 12|3|4\\
R_{2}: & 1|22|3, 1|22|4, 1|23|4\\
R_{111}:& 1|2|3|4
\end{align}
The sum of the coefficients forms Sequence A239204,
$$
1, 1, 4, 17, 76, 353, 1688, 8257, 41128, 207905,\dots
$$
whose generating series is
\begin{equation}\label{eq:rub}
 1+\frac{(x-1)\sqrt{x^2-6 x+1}-x^2-4x+1}{8x^2}. 
\end{equation}

Indeed, 
each $S^I$ is the sum of $2^{\ell(I)-1}$ ribbons $R_J$, so that
the generating series for the sum $a_n$ of the coefficients of $g_n$ on the $R_I$
is obtained from the specialization $S_n\mapsto ux^n$ for $n\ge 1$. This results into the functional equation
\begin{equation}\label{eq:grub}
g = 1+\frac{uxg}{1-xg} = \sum_{n\ge 0}g_n(u)x^n
\end{equation}
and 
\begin{equation}\label{eq:usub}
a_n = \left.\frac{g_n(u)}{u}\right|_{u=2}.
\end{equation}
Solving the quadratic equation \eqref{eq:grub} and making the substitution \eqref{eq:usub}
in the specialization of $(g-1)/(\sigma_x-1)$,
with $\sigma_x=1+ux/(1-x)$,
 we obtain \eqref{eq:rub}.

The coefficients of $g$ in the basis $\Lambda^I$ are, up to sign and
conjugation, the same as in the basis $R_I$:
\begin{equation}
[\Lambda^I]g = (-1)^{|I|-\ell(I)}[R_{I^{\sim}}]g.
\end{equation}
Thus, for example, in
\begin{equation}
\gamma_3= 3\Lambda^3 -6\Lambda^{21} -5\Lambda^{12} + 9\Lambda^{111},
\end{equation}
the coefficients are given by
\begin{align}
-\Lambda^4+4\Lambda^{31}&\rightarrow 3\Lambda^3\\
3\Lambda^{22}-9\Lambda^{211}&\rightarrow -6\Lambda^{21}\\
2\Lambda^{13}-7\Lambda^{121}&\rightarrow -5\Lambda^{12}\\
-3\Lambda^{112}+14\Lambda^{1111}&\rightarrow 9\Lambda^{111}
\end{align}
In the second line, the coefficient of $\Lambda^{21}$ is the difference
between the number of parking quasi-ribbons of shape 121 (3) and those of
shape 31 (9), thus $-6$.

The sum of the absolute values of the coefficients forms (up to a shift) sequence A238112,
$$
1, 1, 5, 23, 107, 509, 2473, 12235, 61463, 312761, 1609005,\ldots
$$
whose generating series is
\begin{equation}
1+\frac{1-5x+2x^2+(2x-1)\sqrt{x^2-6x+1}}{4x^2}.
\end{equation}
This series is obtained by the specialization $S^I\mapsto 2^{|I|-\ell(I)}x^{|I|}$, which amouts to
replacing $g$ by the solution of
\begin{equation}
g(x) = 1+\frac12\sum_{n\ge 1}(2xg)^n = 1+\frac{xg}{1-2xg}
\end{equation}
and $\sigma_x$ by $x/(1-2x)$ in $\gamma=(g-1)/(\sigma_x-1)$.

\section{Gessel's formulas}

Gessel \cite{Ges} has given the formula
\begin{equation}
\gamma=\left(1-\sum_{n\ge 1}S_n\left(1+g+\cdots+g^{n-1}\right)\right)^{-1}
\end{equation}
that is still valid for the noncommutative version. The proof can be copied
verbatim, only replacing $g$ by $\sum S_ng^n$ in this order. His second
formula (for $H$) shows that there is an analog of the geode factorisation for
the generating series $h=1-1/g$ of prime parking functions. This can also be
seen as follows: the operator $S_1^{-1}$ satisfies
\begin{equation}
(uv)S_1^{-1} = u\cdot vS_1^{-1} + uS_1^{-1}v_0
\end{equation}
where $v_0$ is the constant term of $v$. Thus, setting $\eta:=hS_1^{-1}$,
\begin{equation}
\gamma=gS_1^{-1} =(1-h)^{-1}S_1^{-1} = (1-h)^{-1}\cdot hS_1^{-1} =g\eta,
\end{equation}
so that $\eta=g^{-1}\gamma$, and 
\begin{equation}
\eta(\sigma_1-1) = g^{-1}\gamma(\sigma_1-1)=g^{-1}(g-1)=1-g^{-1}=h.
\end{equation}
Of course, by definition, $\eta$ has nonnegative integer coefficients, whose
combinatorial interpretation is immediate, since prime parking functions
correspond to plane trees whose rightmost subtree of the root is a leaf.

\section{The $k$-geodes and beyond}
The $k$-Lagrange series \cite{NTLag2}, defined by
\begin{equation}\label{eq:klag}
g^{(k)} = \sum_{n\ge 0}S_n\left(g^{(k)}\right)^{kn}
\end{equation}
or equivalently by $g^{(k)}=\phi_k(g)$, satisfies a similar property, where
$\phi_k$ is the algebra morphism sending $S_n$ to $S_{n/k}$ if $n$ is
divisible by $k$, and to 0 otherwise. Thus,
\begin{equation}
g^{(k)} = \phi_k(1+\gamma(\sigma_1-1))=1+\phi_k(\gamma)(\sigma_1-1).
\end{equation} 
Note that \cite[Lemma 2.1]{NTLag2}:
\begin{equation}
g^{(k)}=\sum_{n\ge 0}g_n^{(k-1)}\left[g^{(k)}\right]^n
\end{equation}
implies also that $g^{(k)}-1$ is divisible
by $g^{(k-1)}-1$ and that the quotient has nonnegative integer coefficients.

The hierarchy of $k$-Lagrange series can be extended to negative values of $k$, and in fact to any complex value of $k$, since \eqref{eq:klag} still makes sense in this case. Then, $g^{(0)}$ is just
$\sigma_1$, and $g^{(-1)}$ is the series of noncommutative free cumulants introduced in \cite{JVNT}.  
Indeed, by definition,
\begin{equation}
g^{(-1)} = \sum_{n\ge 0}S_n\left[g^{(-1)}\right]^{-n},
\end{equation}
whilst the series $K$ of free cumulants is defined by
\begin{equation}
\sigma_1 = \sum_{n\ge 0}K_n \sigma_1^n,
\end{equation}
which can be interpreting as defining  $g^{(0)}$ in terms of $g^{(-1)}$. It is proved in
\cite[Th. 3.1]{JVNT} that $K=g(-A)^{-1}$, and in \cite[Th. 6.1]{NTLag} that
\begin{equation}
g(-A)^{-1}=\sum_{n\ge 0} S_n g(-A)^n = \sum_{\ge 0}S_n \left[g(-A)^{-1}\right]^{-n},
\end{equation}
so that indeed, $g^{(-1)}=K=g(-A)^{-1}$.

This suggests to introduce the \emph{Lagrange transform} as the algebra
automorphism defined by
\begin{equation}
L(S_n) =g_n, 
\end{equation}
so that $g^{(k)}=L^k(g^{(0)})$ ($k\in\Z$), and to introduce two hierarchies 
$\gamma^{(k)}$ and $\theta^{(k)}$ of higher order geodes, such that
\begin{equation}
g^{(k)}= 1+ \gamma^{(k)}(\sigma_1-1),\quad \gamma^{(k)}=\phi_k(\gamma)
\end{equation}
and
\begin{equation}
g^{(k)}= 1+\theta_k(g^{(k-1)}-1),\quad \theta_k=L^{k-1}(\gamma). 
\end{equation}

The coefficients of $g^{(k)}$ are given by \cite[Eq. (44)]{NTLag}. Indeed, the coefficient of
$S^I$ in $g^{(k)}$ is equal to the coefficient of $S^{kI}=S^{(ki_1,\ldots,ki_p)}$ in $g$, therefore to
\begin{equation}\label{eq:delta}
\delta_I^{(k)}= \sum_{a=(a_1,\ldots,a_{p})}\prod_{j=1}^{p-1}{ki_j\choose a_j}
\end{equation}
where $a$ runs over the Polish codes of all plane tres with $p$ nodes (and therefore ends with a 0
which can be omitted in the product).

For example, writing ${kI\choose a}$ for the product of binomial coefficients
in \eqref{eq:delta}, the coefficient of $S^{211}$ in $g^{(k)}$ is given by
\begin{equation}
{2k\,k\, k\choose 1\,1\,0}+{2k\,k\, k\choose 2\,0\,0} = 2k^2+\frac{2k(2k-1)}{2}=4k^2-k.
\end{equation}
Being polynomials in $t$, the coefficients of $g^{(t)}$ are given by the same formula 
for arbitrary $t$.

\section{Schr\"oder trees and the $e$-Lagrange series}

The binomial coefficients ${t\choose k}$ are the elementary functions $e_k(t)$ when $t$ is a binomial element in the sense of $\lambda$-rings (i.e., a scalar). Thus, if we define a series $g^{[e]}$ by
the functional equation
\begin{equation}\label{eq:defge}
g^{[e]}=\sum_{n\ge 0}S_n\left(\sum_{k\ge 0}e_k \left(g^{[e]}\right)^k\right)^n
\end{equation}
where $e_n$ is an arbitrary sequence (with $e_0=1$) regarded as elementary
functions of a virtual alphabet $\AA$, the coefficient of $S^I$ in $g^{[e]}$
is
\begin{equation}\label{eq:deltage}
\delta_I^{[e]}= \sum_{a=(a_1,\ldots,a_{p})}\prod_{j=1}^{p-1}e_{a_j}(i_j\AA).
\end{equation}
The elementary functions $e_k(m\AA)$ can be expanded as
\begin{equation}
e_k(j\AA) = \sum_{i_1+\cdots+i_j=k}e_{i_1}(\AA)\cdots e_{i_j}(\AA).
\end{equation}
The coefficient of $S^I$ becomes then a sum of products $e_\lambda(\AA)$. By setting $S_n=z^n$
and $e_k=q^k$, \eqref{eq:defge} can be solved explicitely as 
\begin{equation}
g^{[e]} = \frac{1-z-\sqrt{1-2(1+2q)z+z^2}}{2q}
\end{equation}
and one can see that the total number of terms in $g^{[e]}_n$ is the (large) 
Schr\"oder number $s_{n-1}$ (A006318),
and moreover that the sum of these coefficients counts Schr\"oder paths according
to the number of up steps $U=(1,0)$ (A088617). 

The first terms of $g^{[e]}$ are
\begin{align}
g^{[e]}_0 &= 1  \\
g^{[e]}_1 &=  S_1 \\
g^{[e]}_2 &= S_2+e_1 S^{11}  \\
g^{[e]}_3 &= S_3+ 2e_1S^{21}+e_1S^{12}+(e_2+e_{11})S^{111}  \\
g^{[e]}_4 &= S_4+3e_1S^{31}+2e_1S^{22}+(3e_{11}+2e_2)S^{211} \nonumber \\
&+e_1S^{13}+(2e_{11}+e_2)S^{121}+(e_{11}+e_2)S^{112}+(e_{111}+3e_{21}+e_3)S^{1111}
\end{align}

A Schr\"oder tree of size $n$ is a tree with $n+1$ leaves
and whose internal nodes are of arity strictly greater than $1$. A prime
Schr\"oder tree is a Schr\"oder tree whose rightmost subtree of the root is a leaf.
Equivalently, Schr\"oder trees of size $n$ are rooted plane trees with internal
nodes labelled by a composition $I$ of $n$ (read in prefix order) such that a
node labelled  $i$ is of arity $i+1$. Labelling the leaves by 0 and reading the tree in prefix order,
we obtain a Schr\"oder pseudocomposition $I$. The vector space ${\mathcal S}$ spanned by the symbols
$S^I$ is called the Schr\"oder operad \cite{JVNT}.

Schr\"oder paths are in bijection with prime Schr\"oder trees, which are the natural objects for the combinatorial description of $g^{(-1)}$ (free cumulants). In order to recover a similar description for $g^{[e]}$, 
we can lift \eqref{eq:defge}  to the system  
\begin{equation}
	G = (1+X)S_0,\qquad 
\begin{cases}X=\sum_{n\ge 1}S_nY^n\\Y=S_0+\sum_{n\ge 1}e_nX^nS_0.\end{cases}
\end{equation}
Then, for $S_0=1$, $G$ becomes $g^{[e]}$.  Solving recursively, 
we obtain
\begin{align}
Y_0&=E_0=S_0\\
X_1&=S^{10}\\
Y_1 &=  e_1S^{100}\\
X_2&= e_1S^{1100}+S^{200}\\
Y_2 &= e_1^2S^{11000}+e_1S^{20000}+e_2S^{10100}\\
X_3&=e_1^2S^{111000}+e_1S^{12000} \\
&+e_2S^{110100}+e_1S^{21000}+e_1S^{20100}+S^{3000}. \nonumber
\end{align}
Then, it follows by induction that 
\begin{equation}
Y_n = \sum_{t\in \ST_n} e_{\lambda(t)} S^t
\end{equation}
where $\ST_n$ is the set of Schr\"oder trees with $n+1$ leaves, and $\lambda(t)$ is the partition formed by the lengths of the right branches of $t$, where
as in \cite{JVNT}, we set $S^t=S^I$ where $I$ is the code of $t$, and that
\begin{equation}
G_n = 
\sum_{t\in \PST_n}\prod_i e_{\lambda(t_i)} S^t
\end{equation}
where $\PST_n$ denotes the set of prime Schr\"oder trees of size $n$, and the $t_i$
are the subtrees of the root of $t$.

The formula of \cite{JVNT} for free cumulants corresponds to the choice $e_n=(-1)^n$. Indeed,
the coefficient of $S^t$ becomes then $(-1)^{i(t)-1}$, where $i(t)$ is the number of internal nodes of $t$.

Setting $t=1$, we recover the expression of the Lagrange series in terms of Schröder tres given in \cite{MNT}: it is the sum of all trees such that the rightmost subtree of each internal vertex is a leaf. 

For example, the prime Schr\"oder trees of size $3$ are 
\begin{equation}
\begin{split}
&
{\xymatrix@C=2mm@R=2mm{
 *{} & {3}\ar@{-}[dl]\ar@{-}[d]\ar@{-}[dr]\ar@{-}[drrr] \\
 {.} & {.} & {.} & *{} & {.}
}}
\qquad
{\xymatrix@C=2mm@R=2mm{
 *{} & *{} & *{} & {2}\ar@{-}[dl]\ar@{-}[d]\ar@{-}[dr] \\
 *{} & *{} & {1}\ar@{-}[dl]\ar@{-}[dr] & {.} & . \\
 *{} & {.} & *{} & . 
      }}
\qquad
{\xymatrix@C=2mm@R=2mm{
 *{} & *{} & *{} & {2}\ar@{-}[dl]\ar@{-}[d]\ar@{-}[dr] \\
 *{} & *{} & {.} & {1}\ar@{-}[dl]\ar@{-}[dr] & {.} & \\
 *{} & {} & {.} & {} & {.} 
      }}
\qquad
{\xymatrix@C=2mm@R=2mm{
 *{} & *{} & *{} & {1}\ar@{-}[dl]\ar@{-}[dr] \\
 *{} & *{} & {2}\ar@{-}[dl]\ar@{-}[d]\ar@{-}[dr] & {} & {.} & \\
 *{} & {.} & {.} & {.} & { } 
      }}
\\
&
{\xymatrix@C=2mm@R=2mm{
 *{} & *{} & *{} & {1}\ar@{-}[dl]\ar@{-}[dr] \\
 *{} & *{} & {1}\ar@{-}[dl]\ar@{-}[dr] & {} & {.} & \\
 *{} & {1}\ar@{-}[dl]\ar@{-}[dr] & {} & {.} & { }  \\
 {.} & {} & .
      }}
\qquad
{\xymatrix@C=2mm@R=2mm{
 *{} & *{} & *{} & {1}\ar@{-}[dl]\ar@{-}[dr] \\
 *{} & *{} & {1}\ar@{-}[dl]\ar@{-}[dr] & {} & {.} & \\
 *{} & {.} & {} & {1}\ar@{-}[dl]\ar@{-}[dr] & {} &  \\
  {} & {} & . & & .
      }}
\end{split}
\end{equation}
so that
\begin{equation}
G_3=S^{30000}+e_1S^{210000}+e_1S^{201000}+e_1S^{120000}+e_1^2S^{1110000}+e_2S^{1101000}.
\end{equation}

Now, as previously, the multiset of trees obtained by cutting in $G_{n+k}$ the rightmost $k$-corolla
followed only by leaves of the roots when visited in prefix order does not
depend of $k$, and as the coefficient of a tree depends only of the lengths of
the right branches, the $e$-geode is well defined, and given by
\begin{equation}
\gamma^{[e]}= g^{[e]}S_k^{-1}
\end{equation}
for any $k\ge 1$. For example,
\begin{align}
\gamma_1^{[e]}&=e_1S_1\\
\gamma_2^{[e]}&=2e_1S^2+(e_1^2+e_2)S^{11}\\
\gamma_3^{[e]}&= 3e_1S^3+(3e_{1}^2+2e_2)S^{21}+(2e_{1}^2+e_2)S^{12}+(e_{1}^3+3e_{1}e_{2}+e_3)S^{111}\\
\gamma_4^{[e]}&=
 4 e_{1} S_{4} 
+ \left(6 e_{1}^{2} + 3 e_{2}\right) S^{31} 
+ \left(5 e_{1}^{2} + 2 e_{2}\right) S^{22}\\
&+ \left(4 e_{1}^{3} + 8 e_{1} e_{2} + 2 e_{3}\right) S^{211}
+ \left(3 e_{1}^{2} + e_{2}\right) S^{13} \\
&+ \left(3 e_{1}^{3} + 5 e_{1} e_{2} + e_{3}\right) S^{121} 
+ \left(2 e_{1}^{3} + 4 e_{1} e_{2} + e_{3}\right) S^{112}\\
&+ \left(e_{1}^{4} + 6 e_{1}^{2} e_{2} + 4 e_{1} e_{3} + 2 e_{2}^{2} + e_{4}\right) S^{1111}
\end{align}

\section{Tables}

Here are the expressions of $g_k^{(t)}$, $\gamma_k^{(t)}$,
$\theta_k^{(t)}$, $h_k^{(t)}$, and $\eta_k^{(t)}$ for $k\leq4$.
{\footnotesize
\begin{align}
g^{(t)}_1 &=S_1\\
g^{(t)}_2 &=S_2+tS^{11}\\
g^{(t)}_3 &=S_3+2tS^{21}+tS^{12}+\frac{3t^2-t}{2}S^{111}\\
g^{(t)}_4 &=S_4+3tS^{31}+2tS^{22}+(4t^2-t)S^{211}+tS^{13}+\frac{5t^2-t}{2}S^{121}\\
&+\frac{3t^2-t}{2}S^{112}+\frac{8t^3-6t^2+t}{3}S^{1111}\nonumber
\end{align}

\begin{align}
\gamma^{(t)}_1 &=tS_1\\ 
\gamma^{(t)}_2 &= 2tS_2+\frac{3t^2-t}{2}S^{11}\\
\gamma^{(t)}_3 &= 3tS_3+(4t^2-t)S^{21}+\frac{5t^2-t}{2}S^{12}+\frac{8t^3-6t^2+t}{3}S^{111}\\ 
\gamma^{(t)}_4 &= 4tS_4
+\frac{15t^2-3t}{2}S^{31} +(6t^2-t)S^{22}+\frac{25t^3-15t^2+2t}{3}S^{211} +\frac{7t^2-t}{2}S^{13}\\
&+\frac{17t^3-9t^2+t}{3}S^{121} +\frac{25t^3-15t^2+2t}{6}S^{112}\nonumber\\
&+\frac{125t^4-150t^3+55t^2-6t}{24}S^{1111}\nonumber
\end{align}

\begin{align}
\theta_1^{(t)} &= S_1\\
\theta_2^{(t)} &= 2S_2+(2t-1)S^{11}\\ 
\theta_3^{(t)} &= 3S_3+(6t-3)S^{21}+(3t-1)S^{12}+\frac{9t^2-11t+4}{2}S^{111}\\ 
\theta_4^{(t)}&= 4S_4+(12t-6)S^{31}+(8t-3)S^{22}+(16t^2-19t+7)S^{211}+(4t-1)S^{13}\\
&+(10t^2-10t+3)S^{121}+(6t^2-6t+12)S^{112}\nonumber\\
&+\frac{64t^3-129t^2+101t-30}{6}S^{1111}\nonumber
\end{align}

\begin{align}
h^{(t)}_1 &= S_1\\
h^{(t)}_2 &= S_2\\
h^{(t)}_3 &= S_3+tS^{21}\\
h^{(t)}_4 &=S_4+2tS^{31}+tS^{22}+\frac{3t^2-t}{2}S^{111}
\end{align}

\begin{align}
\eta^{(t)}_1 &= 0\\
\eta^{(t)}_2 &=tS_2\\
\eta^{(t)}_3 &=2tS_3+\frac{3t^2-t}{2}S^{21}\\
\eta^{(t)}_4 &=3tS_4+(4t^2-1)S^{31}+\frac{5t^2-t}{2}S^{22}+\frac{8t^3-6t^2+t}{3}S^{1111}
\end{align}
}


\footnotesize
\begin{thebibliography}{aa}
%
\bibitem{EH}{ R. Ehrenborg} and {A. Happ},
{ The antipode of the noncrossing partition lattice},
Advances in Applied Mathematics,
110 (2019),  76--85.
%
\bibitem{Einz}{ H. Einziger},
{ Incidence Hopf Algebras: Antipodes, Forest Formulas, and Noncrossing
Partitions},
Thesis (Ph.D.) -- The George Washington University, ProQuest LLC, Ann Arbor, MI, 2010. 111 pp.
ISBN:978-1124-15739-9.



%
\bibitem{NCSF1}{ I. M. Gelfand, D. Krob, A. Lascoux, B. Leclerc,
V.~S. Retakh}, and { J.-Y. Thibon},
{ Noncommutative symmetric functions},
Adv. in Math. 112 (1995), 218--348.
%
\bibitem{Ges} I. M. Gessel, Lattice paths and the geode, arXiv:2507.09405.

\bibitem{Gos} F. Gossow, Ordered trees and the geode, arXiv:2507.18097.
%
\bibitem{GJ1}{ I. P. Goulden} and { D. M. Jackson},
{ The Combinatorial Relationship Between Trees, Cacti and Certain
Connection Coefficients for the Symmetric Group},
Europ. J. Combinatorics { 13} (1992), 357--365.
%
\bibitem{GJ2}{ I. P. Goulden} and { D. M. Jackson},
{ Symmetric functions and Macdonald's result for top connexion coefficients in
the symmetric group},
J. Algebra { 166} (1994), no. 2, 364--378.
%
\bibitem{Hai1}{ M. Haiman},
{ Conjectures on the quotient ring by diagonal invariants},
J. Algebraic Combin.  3 (1994), 17--36.
%
\bibitem{HNTtrees}{F. Hivert, J.-C. Novelli} and { J.-Y. Thibon},
Trees, functional equations, and combinatorial Hopf algebras,
European J. Combin. 29 (2008), no. 7, 1682--1695.
%
\bibitem{JVNT} M. Josuat-Vergès, F. Menous, J.-C. Novelli, and J.-Y. Thibon, 
Free cumulants, Schröder trees, and operads, 
Advances in Applied Mathematics 88 (2017), 92--119.
%
%
\bibitem{Mcd}{ I. G. Macdonald},
{ Symmetric functions and Hall polynomials},
2nd ed., Oxford University Press, 1995.
%
\bibitem{MNT}
F. Menous, J.-C. Novelli and J.-Y. Thibon, Combinatorics of Poincar\'e's and Schr\"oder's equations, in {\it Resurgence, physics and numbers}, 329--378, CRM Series, 20, Ed. Norm., Pisa.


\bibitem{NTLag}{ J.-C Novelli} and { J.-Y. Thibon},
{ Noncommutative symmetric functions and Lagrange inversion},
Adv. Appl. Math. { 40} (2008), 8--35.
%
\bibitem{NTLag2}{ J.-C Novelli} and { J.-Y. Thibon},
{ Noncommutative symmetric functions and Lagrange inversion II:
noncrossing partitions and the Farahat-Higman algebra},
Adv. in Appl. Math. 140 (2022), Paper No. 102396, 39 pp..

%
\bibitem{NTDup} J.-C. Novelli and J.-Y. Thibon,
 Duplicial algebras, parking functions, and Lagrange inversion, in {\it Algebraic combinatorics, resurgence, moulds and applications (CARMA). Vol. 1}, 263--290, IRMA Lect. Math. Theor. Phys., 31, EMS Publ. House, Berlin.
%
\bibitem{NTpark}{ J.-C. Novelli}, { J.-Y. Thibon},
{ Hopf algebras and dendriform structures arising from parking functions},
Fund. Math.   193  (2007),   189--241.
%
\bibitem{Ran}{ G. N. Raney}, 
{ Functional composition patterns and power series reversion},
Trans. Amer. Math. Soc.  94 (1960), 441--451.
%
%
\bibitem{Slo}{ N. J. A. Sloane},
{ The On-Line Encyclopedia of Integer Sequences},\\
http://www.research.att.com/~njas/sequences/.
%
\bibitem{WR}
 N. J. Wildberger and Dean Rubine, A hyper-Catalan series solution to polynomial equations, and the
Geode, Amer. Math. Monthly 132 (2025), no. 5, 383--402.
%
%
%
\end{thebibliography}
\end{document}